\documentclass[11pt,leqno]{article}
\usepackage{amssymb}
\usepackage{graphicx}
\usepackage{epic,eepic}
\newcommand{\ds}{\displaystyle}

\newcommand{\Z}{\mathbb{Z}}

\newcommand{\ol}{\overline}

\newcommand{\ra}{\rightarrow}

\newcommand{\name}[2]{{#1}{\scriptsize{#2}}}
\renewcommand{\name}[1]{{\sc #1}}

\begin{document}

\begin{center}
{\Large \bf A catalogue of small regular matroids and their Tutte polynomials

}

\medskip
\name{Harald} \name{Fripertinger} \name{and} \name{Marcel} \name{Wild}
\end{center}

\begin{abstract}\noindent
A catalogue of all non-isomorphic simple connected regular matroids ${\cal M}$ of cardinality $n \leq 15$ is provided on the net. These matroids are given as binary matrix matroids and are sieved from the large pool of all non-isomorphic binary matrix matroids of cardinality $\leq 15$.  For each ${\cal M}$ its Tutte polynomial is determined by an algorithm based on internal and external base activity. 
\end{abstract}

\section{Introduction}

We assume familiarity with basic matroid concepts such as contraction
matroids, dual matroids, or the flat lattice \cite{[5]}. A famous binary matroid is the {\it Fano} matroid $F_7$ on the set $\{a,b, \cdots, g\}$ which can be defined as the column matroid of this $3 \times 7$ matrix over $GF(2)$:

\begin{center}
$A: =$ \quad \begin{tabular}{|c|c|c|c|c|c|c|c}
$a$  & $b$ & $c$ & $d$ & $e$ & $f$ & $g$ & \\ \hline
$1$ & $0$ & $0$ & $0$ & $1$ & $1$ & $1$ & $=r_1$\\ \hline
$0$ & $1$ & $0$ & $1$ & $0$ & $1$ & $1$ & $=r_2$\\ \hline
$0$ & $0$ & $1$ & $1$ & $1$ & $0$ & $1$ & $=r_3$ \\ \hline \end{tabular}
\end{center}

It is essential in the sequel that two binary matroids are isomorphic if and
only if some obviously sufficient condition holds for any two respective
matrix representations. Specifically, for any column representation of a
binary $n$-matroid, {\it every} other column representation is obtained by
replacing the rows $r_i$ of the matrix by a set of equivalent rows, i.e.
generating the same subspace of $GF(2)^n$, and by subsequently permuting the
labelled columns \cite[10.1.3]{[5]}. Thus, in our case we could take this representing matrix for $F_7$:

\begin{center}
\begin{tabular}{|c|c|c|c|c|c|c|l}
$a$ & $b$ & $c$ & $d$ & $e$  & $f$ & $g$ & \\ \hline
$0$ & $1$ & $0$ & $1$ & $0$ & $1$ & $1$  & $=r_2$\\ \hline
$1$ & $0$ & $1$ & $1$ & $0$ & $1$ & $0$ & $=r_1+r_3$\\ \hline
$0$ & $1$ & $1$ & $0$ & $1$ & $1$ & $0$ & $=r_2+r_3$\\ \hline \end{tabular}
\end{center}
The fact that this, like every matrix representation of $F_7$,  has the same column set as $A$ is particular to $F_7$. It makes it easy to decide whether or not a $3 \times 7$-matrix over $GF(2)$ represents~$F_7$.  

It is well known that the dual ${\cal M}^d$ (here $= F^d_7$)   of a representable matroid ${\cal M}$  is obtained by taking the orthogonal complement $Y = X^\perp$ of $X =$ rowspace$(A)$ with respect to the canonical scalar product $r \cdot s$ in $GF(2)^7$, and letting $A'$ be any matrix whose rows constitute a basis of $Y$.  For instance, one verifies that $r_i \cdot s_j =0$ for all $1 \leq i \leq 3$ and $1 \leq j \leq4$:

\begin{center}
$A' =$ \quad \begin{tabular}{|c|c|c|c|c|c|c|c} 
$a$ & $b$ & $c$ & $d$ & $e$ & $f$ & $g$ & \\ \hline
$0$ & $1$ & $1$ & $1$ & $0$ & $0$ & $0$ & $=s_1$ \\ \hline
$0$ & $1$ & $1$ & $0$ & $1$ & $1$ & $0$ & $=s_2$ \\ \hline
$1$ & $1$ & $0$ & $0$ & $0$ & $1$ & $0$ & $=s_3$ \\ \hline
$1$ & $1$ & $1$ & $0$ & $0$ & $0$ & $1$ & $=s_4$ \\ \hline  \end{tabular}
\end{center}

Therefore, the column matroid of $A'$ must be isomorphic to $F^d_7$. In fact, whether or not {\it any} $4 \times 7$ matrix of rank 4, with columns labelled $a, b, \ldots, g$ from left to right, yields a column matroid isomorphic to $F^d_7$, is most easily checked by testing whether its rows are orthogonal to $r_1, r_2, r_3$. 

The {\it simplification} $\ol{\cal M}$ of a matroid ${\cal M}$ is the matroid formed by removing all loops and identifying any multiple points into a single point. Recall that a matroid is {\it regular} if it is representable over every field. There are various ways to characterize the regular matroids among the class of binary matroids. For us the algorithmic most convenient one was the following.

\medskip
{\bf Theorem 1 \cite[13.4.1(ii)]{[5]}}: {\em  A binary matroid ${\cal M}$ of
rank $k$ is regular if and only if there neither is a rank $k-3$ flat $U$ with
$\ol{{\cal M} / U} \simeq F_7$, nor a rank $k-4$ flat $V$ with $\ol{{\cal M} /
V} \simeq F_7^d$.}

\medskip
In view of Theorem 1 the following procedure suggests itself to generate up to isomorphism all regular matroids up to size $N$:
\begin{enumerate}
 \item [a)] Generate up to isomorphism all binary matroids ${\cal M}$ up to size $N$.
\item[b)] For each fixed ${\cal M}$ of rank $k \leq N$ compute all rank $k-3$ flats $U_1, \cdots, U_s$ and all rank $k-4$ flats $V_1, \cdots, V_t$.
\item[c)] Check whether $\ol{{\cal M}/U_i} \simeq F_7$ for some $1 \leq i \leq s$, or $\ol{{\cal M}/V_j} \simeq F^d_7$ for some $1 \leq j \leq t$.
\end{enumerate}

Here comes the section break up. The details of a), b), c) are discussed in sections 2, 3, 4 respectively. Section 5 outlines a novel algorithm for computing the Tutte polynomial of a matroid. It was used to calculate the Tutte-polynomials for all regular matroids generated.

\section{Generating binary matroids up to isomorphism}
Expanding upon the introductory remarks, each binary matrix matroid ${\cal M}$
of cardinality $n$ and rank $k$ can be represented by a binary 
$k\times n$-matrix $A$ of rank $k$. This representation, however, is not
unique, since all matrices of the form $G\cdot A$ for $G\in {\rm GL}_k(2)$, the
group of all invertible $k\times k$ matrices over $GF(2)$, represent the same
matroid. 

Two matrix matroids represented by the $k\times n$-matrices $A$ and
$B$ over $GF(2)$ of rank $k$ are isomorphic 
if and only if there exists an $n\times
n$-permutation matrix $P$
and some $G\in {\rm GL}_k(2)$ so that 
$B=G\cdot A\cdot P$. 
Thus, the isomorphism class of the matrix matroid given 
by $A$ is the orbit of $A$ under the action of ${\rm GL}_k(2)\times S_n$ on
the set of all $k\times n$-matrices of rank $k$. This orbit consists of all
matrices of the form $G\cdot A\cdot P$
where $G\in {\rm GL}_k(2)$ and $P$ is a 
permutation matrix.

Next we consider these matrices $A$ as functions from ${\bf n}: = \{1, \ldots,
n\}$ to $GF(2)^k$, where $A(i)^\top$,
the transposed of $A(i)$, is the $i$-th column
of the matrix $A$. The multiplication $A\cdot P$
of matrices is then
replaced by composing $A$ with a suitable permutation $\pi\in S_n$. 
Thus we have a group action of ${\rm GL}_k(2)\times S_n$ on the set of all
functions from ${\bf n}$ 
to $GF(2)^k$ determining matrices of rank $k$.
Since the symmetric group $S_n$ acts on the domain of these functions, it is
not important in which position a vector  $v\in GF(2)^k$ appears in $A$.  We
are only interested in how often $v$ appears in $A$. 
The group action of ${\rm GL}_k(2)\times S_n$ 
can be replaced by an action of ${\rm GL}_k(2)$ on the set of
all $S_n$-orbits of functions mapping ${\bf n}$
to $GF(2)^k$.
The $S_n$-orbit of a function (or matrix) $A$ is a multiset of
elements of $GF(2)^k$ of size $n$. It contains each column of the matrix $A$
together with its multiplicity.
Hence, for obtaining the
isomorphism classes of binary matroids it is enough 
to consider the action of ${\rm GL}_k(2)$ on the set of all multisets of
elements of $GF(2)^k$ of size $n$. Such multisets are written as functions 
$f: GF(2)^k \to \Z_{\geq 0}:=\{0,1,2, \ldots \}$
where $f(v)$ is the multiplicity of the
column $v\in GF(2)^k$ in $f$. Therefore, we call $f$ the multiplicity function
of a multiset of $GF(2)^k$ of size $n$. Consequently, the multiplicity
functions $f$ of a multiset of $GF(2)^k$ of rank $k$ and size $n$
satisfy the following two properties:
\begin{enumerate} 
\item [a)] There exist $k$ linearly independent vectors 
$v \in GF(2)^k$ such that $f(v) > 0$ and 
\item[b)] $\sum_{v\in GF(2)^k} f (v)=n$.
\end{enumerate}

Finally, we introduce the function $\rho: 
GF(2)^k\to \{0,1,\ldots,2^k-1\}$ for
labelling the elements $v=(v_1,\ldots,v_k)\in GF(2)^k$ with nonnegative
integers. (Thus $\rho$ is a rank function on $GF(2)^k$.) The label of $v$ is
$$\rho(v_1,\ldots,v_k):=\sum_{j=1}^kv_j2^{j-1},$$
where we identify the elements $0$ and $1$ of $GF(2)$ with the corresponding
integers. 
The unit vectors
$(1,0,\ldots,0)$, $(0,1,0,\ldots,0)$, \ldots, 
$(0,\ldots,0,1)$ in $GF(2)^k$
are then labelled by the powers of $2$, i.e. by
$1$, $2$, \ldots, $2^{k-1}$, respectively. 

Conversely, given a label $r\in\{0,1,\ldots,2^k-1\}$ we determine the unrank
function $\rho^{-1}:\{0,1,\ldots,2^k-1\}\to GF(2)^k$ by
$\rho^{-1}(r)=v=(v_1,\ldots,v_k)$ which is obtained from the binary
representation of $r$ as $r=\sum_{j=1}^kv_j2^{j-1}$ with $v_j\in\{0,1\}$.

Using $\rho$, we determine a permutation representation $\tilde {\rm GL}_k(2)$
of ${\rm GL}_k(2)$ on the set $\{0,1,\ldots,\linebreak[3]2^k-1\}$. 
Consider $B\in {\rm
GL}_k(2)$ and $j\in \{0,1,\ldots,2^k-1\}$, then $\pi_B$ defined by
$\pi_B(j):={\rho(B\cdot (\rho^{-1}(j))^\top )}$ is a permutation of
$\{0,1,\ldots,2^k-1\}$.

Using $\rho$ we consider the multiplicity
functions $f$ of multisets of $GF(2)^k$ of size $n$
as functions $f$ from $\{0,1,\ldots,2^k-1\}$ to $\Z_{\geq 0}$ with the 
property $\sum_{j=0}^{2^k-1} f(j)=n$.
We write these functions $f$ as vectors
$f=(f(0),f(1),\ldots,f(2^k-1))$.  The set of these vectors is totally ordered by
the lexicographic order.  The group $\tilde {\rm GL}_k(2)$ acts in a natural
way on the domain of these functions.  This induces a group action on the set
of all functions from $\{0,1,\ldots,2^k-1\}$ to $\Z_{\geq 0}$.  We consider
the lexicographic largest element of the $\tilde {\rm GL}_k(2)$-orbit of $f$
as the standard representative of this orbits. 

The functions $f$ are supposed to describe matrix matroids of rank $k$,
whence, these matrices contain a set of $k$ linearly independent vectors. Due
to the labelling of the vectors of $GF(2)^k$ and since the standard
representative is the largest vector in its orbit, a standard representative
contains the labels of the $k$ unit vectors in $GF(2)^k$ with nonzero
multiplicities.  Thus for $i\in\{2^j\mid 0\leq j<k\}$ we have $f(i)>0$.
Moreover, it is easy to see that standard representatives $f$ satisfy $f(i)\geq
f(r)$ for all $r>i$ and all  $i\in\{2^j\mid 0\leq j<k\}$.

In order to determine a complete list of representatives of binary matrix
matroids of size $n$ and rank $k$ we apply Read's method of orderly generation
(cf. \cite{Read78,ColRead79,ColRead79a,Kerber99}) in the following setting. The group $\tilde {\rm GL}_k(2)$
acts on the domain of all multiplicity functions $f$ from $\{0,1,\ldots,2^k-1\}$ to
$\Z_{\geq 0}$ which satisfy $\sum_{j=0}^{2^k-1} f(j)=n$.  We are listing
all these functions lexicographically decreasing starting with the
lexicographically largest one. Due to our construction we restrict ourselves
to those $f$ with $f(i)>0$ and $f(i)\geq f(r)$ for all $r>i$ and
$i\in\{2^j\mid 0\leq j<k\}$.  Each of these functions is tested whether it is
the lexicographically largest element in its $\tilde {\rm GL}_k(2)$-orbit. If
so, it is considered to be a canonic representative of a binary matrix
matroid. Otherwise it represents a matrix matroid which is isomorphic to a
matrix matroid which was already found before. In this case it is sometimes
possible to obtain information which multiplicity function $f$ should be tested in the next
step (for more details see~\cite{Grund90}).
In several situations we do not proceed with the
lexicographically next $f$ but we are allowed to do some jumps in the set of
all these functions.

A matrix representation $A=A_f$ of a multiplicity function $f$ associated with
a matrix matroid of rank $k$ and size $n$
is then obtained from $f$, by building a $k\times n$-matrix
consisting of $f(i)$ columns of $(\rho^{-1}(i))^\top$ for $i\in
\{0,1,\ldots,2^k-1\}$. It is easier to write down just the labels 
of the columns of this matrix, thus we obtain a vector
$r_f:=(r_1,\ldots,r_n)$ of labels with $r_\nu\leq r_{\nu+1}$ for $\nu<n$ and
$\vert\{\nu\mid 1\leq \nu\leq n,~r_\nu=i\}\vert =f(i)$ for $i\in
\{0,1,\ldots,2^k-1\}$. Due to this construction, if $f_1$ and $f_2$ are
multiplicity functions of size $n$ and $f_1$ is larger than $f_2$ with respect
to the lexicographic order, then the vector $r_{f_1}$ is smaller than
$r_{f_2}$  in the lexicographic order of all vectors in
$\{0,1,\ldots,2^k-1\}^n$.  Since we present the binary matroids in form of
these vectors $(r_1,\ldots,r_n)$, the standard representative of an
isomorphism class is now the smallest element in its orbit and the
representatives are listed monotonically increasing according to the
lexicographic order.

All these computations were done in SYMMETRICA~\cite{Symmetrica}.
For obvious reasons we
restrict ourselves to loopless matroids. Thus the columns of the matrix
matroids belong to $GF(2)^k\setminus\{0\}$, therefore, $f(0)=0$ for all $f$.

This way we obtain complete lists of representatives of the isomorphism
classes of all loopless binary matrix matroids and of
all simple loopless binary matrix matroids of size $\leq 15$ and rank $\leq
7$.

For example, there are two simple loopless binary matrix matroids of 
rank $3$ and size $4$. They are given by 
$$(1,2,3,4) \quad \mbox{which represents the
matrix} \quad \left( \begin{array}{cccc} 1&0&1&0\\ 0&1&1&0\\ 0&0&0&1 
\end{array}\right)$$
and by $$(1,2,4,7) \quad \mbox{which represents the matrix} \quad \left(
\begin{array}{cccc} 1&0&0&1\\ 0&1&0&1\\ 0&0&1&1\end{array} \right).$$
These matroids were found from the list of all binary matroids of rank $3$ and
size $4$. It contains the three vectors 
$r_1:=(1,1,2,4)$, $r_2:=(1,2,3,4)$, $r_3:=(1,2,4,7)$.
The unit vectors have labels 1, 2, 4. Since the rank $k=3$, the multiplicity
functions are functions from $\{0,1,\ldots,7\}$ to $\Z_{\geq 0}$ satisfying
$\sum_{j=0}^7 f(j)=4$. 
The corresponding multiplicity functions are
$f_1=(0,2,1,0,1,0,0,0)$,
$f_2=(0,1,1,1,1,0,0,0)$,
$f_3=(0,1,1,0,1,0,0,1)$ and the corresponding multisets of $GF(2)^3$ of rank
$3$ and size $4$ are 
$A_1=\{(1,0,0),(1,0,0),(0,1,0),(0,0,1)\}$, 
$A_2=\{(1,0,0),(0,1,0),(0,1,1),(0,0,1)\}$, and
$A_3=\{(1,0,0),(0,1,0),(0,0,1),(1,1,1)\}$.

\section{Getting the flats of rank $k-3$ and $k-4$}

Suppose ${\cal M}$ is isomorphic to the column matroid of the $k \times n$ matrix $A$ over $GF(2)$. 
The {\it circuits} of ${\cal M}$ are well known to be the inclusion-minimal
members among the sets $\mbox{supp}(z): = \{i \in {\bf n}\mid
z_i =1\}$ where $z$ ranges over the subspace (rowspace$(A))^\perp$. Actually, circuits come up only in section 5. For the time being we rather need the (say) $h$ many {\it cocircuits} $D_i$ of ${\cal M}$ (i.e. the circuits of ${\cal M}^d$) which can be found as the minimal members among the sets $\mbox{supp}(z)$ where $z$ ranges over the $2^k$-element rowspace of $A$. As for any matroid (representable or not), the complements $H_i: = {\bf n} \setminus D_i \ (1 \leq i \leq h)$ are precisely the hyperplanes of ${\cal M}$, i.e. the rank $k-1$ flats. 

The rank $k-2$ flats $T$ clearly are the inclusionwise maximal sets among the sets $H_i \cap H_j \ (1 \leq i < j \leq h)$. Similarly the rank $k-3$ flats $U$ are obtained from the $T$'s, and then the rank $k-4$ flats $V$ in the same manner from the $U$'s. 


\section{Calculations of contractions}

Let ${\cal M}$ be a binary matroid and $U$ some flat such that the simplification $\ol{{\cal M} / U}$ of the contraction ${\cal M} / U$ is isomorphic to $F_7$. If $I \subseteq U$ is any generating set, i.e. its closure $cl (I)$ is $U$, then also $\ol{{\cal M}/I} = F_7$. It turns out that {\it independent} generating sets $I \subseteq U$ are the most appropriate since matrix representations of ${\cal M}$ smoothly transform into matrix representations of ${\cal M}/ I$.

\medskip\noindent
{\bf Example:} Each binary matroid ${\cal M}$ on the set $\{ a, b, \cdots, m, n\}$, with basis $B = \{a, b, c, d, e \}$, can be represented by a matrix of type

\begin{center}
$(I_5, A)=$ \quad \begin{tabular}{|c|c|c|c|c|c|c|c|c|c|c|c|c|} 
$a$ & $b$ & $c$ & $d$ & $e$ & $f$ & $g$ & $h$ & $i$  & $k$ & $l$ & $m$ & $n$\\ \hline
$1$ & $0$ & $0$ & $0$ & $0$ & $1$ & $0$ & $1$ & $1$ & $0$ & $1$ & $1$ & $0$\\ \hline
$0$ & $1$ & $0$ & $0$ & $0$ & $1$ & $0$ & $1$ & $0$ & $0$ & $0$ & $1$ & $1$ \\ \hline
$0$ & $0$ & $1$ & $0$ & $0$ & $0$ & $1$ & $0$ & $1$ & $0$ & $1$ & $0$ & $0$ \\ \hline
$0$ & $0$ & $0$ & $1$ & $0$ & $0$ & $0$ & $1$ & $1$ & $1$ & $0$ & $1$ & $1$ \\ \hline
$0$ & $0$ & $0$ & $0$ & $1$ & $1$ & $1$ & $1$  & $0$ & $1$ & $1$ & $0$ & $0$ \\ \hline
\end{tabular}
\end{center}
where $I_5$ is the $5 \times 5$ identity matrix. Provided the rows are labelled by the base elements in the right order, this matrix can without loss of information be replaced by the {\it standard matrix} $A$ of ${\cal M}$ with respect to the basis $B$:

\begin{center}
$A=$ \quad \begin{tabular}{c|c|c|c|c|c|c|c|c|}
& $f$ & $g$ & $h$ & $i$ & $k$ & $l$ & $m$ & $n$ \\ \hline
$a$ & $1$ & $0$ & $1$ & $1$ & $0$ & $1$ & $1$ & $0$\\ \hline
$b$ & $1$ & $0$ & $1$ & $0$ & $0$ & $0$ & $1$ & $1$ \\ \hline
$c$ & $0$ & $1$ & $0$ & $1$ & $0$ & $1$ & $0$ & $0$ \\ \hline
$d$ & $0$ & $0$ & $1$  & $1$ & $1$ & $0$ & $1$ & $1$ \\ \hline
$e$ & $1$ & $1$ & $1$ & $0$ & $1$ & $1$ & $0$ &$0$ \\ \hline \end{tabular}
\end{center}
For instance, one still reads off from $A$ that $f = a +b +e$. It turns out that when the independent set $I$ happens to be a subset of $B$, then $B-I$ is a basis of ${\cal M} / I$, and the standard matrix $A^\ast$ of ${\cal M} / I$ with respect to the basis $B-I$ is obtained from $A$ by cancelling the rows with labels in $I$. 

But what if the independent generating set $I'$ of the flat $U$ is not a subset of $B$? Then we extend $I'$ to some other basis $B'$ of ${\cal M}$. To fix ideas, say $I' = \{b, g\} \not\subseteq B$. We switch $g$ with any element in $B$ that features in the column of $g$, say with $c$. By pivoting with respect to $x_{c,g} =1$ the matrix $A$ transforms to some matrix $A'$:

\begin{center}
\begin{tabular}{c|c|c|c|c|c|c|c|c|} 
& $f$ & $g$ & $h$ & $i$ & $k$ & $l$ & $m$ & $n$ \\ \hline
$a$ & $1$ & ${\bf 0}$ & $1$ & $1$ & $0$ & $1$ & $1$ & $0$ \\ \hline 
$b$ & $1$ & ${\bf 0}$ & $1$ & $0$ & $0$ & $0$ & $1$ & $1$ \\ \hline
$c$ & ${\bf 0}$ & ${\bf 1}$ & ${\bf 0}$ & ${\bf 1}$ & ${\bf 0}$ & ${\bf 1}$ & ${\bf 0}$ & ${\bf 0}$ \\ \hline
$d$ & $0$ & ${\bf 0}$ & $1$ & $1$ & $1$ & $0$ & $1$ & $1$ \\ \hline
$e$ & $1$ & ${\bf 1}$ & $1$ & $0$ & $1$ & $1$ & $0$ & $0$ \\ \hline \end{tabular} \quad $\ra$ \quad 
\begin{tabular}{c|c|c|c|c|c|c|c|c|} 
& $f$ & $c$ & $h$ & $i$ & $k$ & $l$ & $m$ & $n$ \\ \hline
$a$ & $1$ & $0$ & $1$ & $1$ & $0$ & $1$ & $1$ & $0$ \\ \hline 
$b$ & $1$ & $0$ & $1$ & $0$ & $0$ & $0$ & $1$ & $1$ \\ \hline
$g$ & $0$ & $1$ & $0$ & $1$ & $0$ & $1$ & $0$ & $0$ \\ \hline
$d$ & $0$ & $0$ & $1$ & $1$ & $1$ & $0$ & $1$ & $1$ \\ \hline
$e$ & $1$ & $1$ & $1$ & $1$ & $1$ & $0$ & $0$ & $0$ \\ \hline \end{tabular}
\end{center}

Generally, the matrix $C' = (y_{\gamma, \delta})$ arising from $C = (x_{\gamma, \delta})$ by {\it pivoting} with respect to $x_{\alpha, \beta} \neq 0$ is defined [Ox, p.84] by
$$y_{\gamma, \delta} : = \left\{ \begin{array}{lll} x_{\gamma, \delta} + x_{\gamma, \beta} x_{\alpha, \delta}, &  \mbox{if} \  \gamma \neq \alpha \ \mbox{and} \ \delta \neq \beta; \\
\\
x_{\gamma, \delta}, & \mbox{otherwise}. \end{array} \right.$$
It turns out that $A'$ above is the standard matrix of ${\cal M}$ with respect to the basis $B': = (B-\{g\}) \cup \{c\} = \{a, b, g, d, e\}$ of ${\cal M}$. Since now $I' \subseteq B'$, the standard matrix $A^\ast$ of ${\cal M} / I'$ with respect to the basis $B'-I' = \{a, d, e\}$ is

\begin{center}
$A^\ast =$ \ \begin{tabular}{c|c|c|c|c|c|c|c|c|}
 & $f$ & $c$ & $h$ & $i$ & $k$ & $l$ & $m$ & $n$  \\ \hline
$a$ & $1$ & $0$ & $1$ & $1$ & $0$ & $1$ & $1$ & $0$ \\ \hline
$d$ & $0$ & $0$ & $1$ & $1$ & $1$ & $0$ & $1$ & $1$ \\ \hline
$e$ & $1$ & $1$ & $1$ & $1$ & $1$ & $0$ & $0$ & $0$ \\ \hline \end{tabular}
\end{center}
From this one reads off that $h, i$ are parallel in ${\cal M} / I'$ and that
$\ol{{\cal M} / I'} \simeq F_7$. Therefore, ${\cal M}$ is not regular by
Theorem 1. \hfill $\square$

\medskip
We mention that it pays off to handle at once all rank $k-3$ flats $U$ of ${\cal M}$ that ``fit'' the current standard matrix $A$ of ${\cal M}$ with respect to base $B$. Namely, they  are the ${k \choose 3}$ many flats $U = cl (I)$, where $I$ ranges over all $(k-3$)-element subsets of $B$. 


\section{Computing the Tutte polynomial}

For an element $a$ of a matroid ${\cal M}$ we denote by ${\cal M} -a$ and ${\cal M}/a$ the deletion and contraction of $a$ respectively. Usually the Tutte polynomial $T({\cal M}, x,y)$ of ${\cal M}$ is calculated via  the deletion-contraction formula
\begin{equation}
T({\cal M}, x, y )\quad =\quad T({\cal M}-a, x,y)\quad + \quad T({\cal M}/a, x,y) 
\label{1}
\end{equation}
which holds for all $a \in {\cal M}$ which neither are isthmuses nor loops. A
sophisticated application of $(1)$ is discussed in \cite{[6]} but it only concerns graphic matroids.

The second author of the present article wrote a very different Mathematica program to calculate $T({\cal M}, x,y)$ for any matroid ${\cal M}$ whose circuits and cocircuits are known. It is based on the formula
\begin{equation}
T({\cal M}, x,y)\quad =\quad \ds\sum_B x^{i (B)} y^{e (B)}
\label{2}
\end{equation}
from \cite[p.~236]{[9]}, 
where the sum is over all bases $B$ of ${\cal M}$ and where $i(B)$, $e(B)$ are the {\it internal} and {\it external} activities of $B$. 

Let us provide the gist of the algorithm; a detailed discussion will be given
elsewhere. First the {\it external activity} of $B$ is defined
\cite[p.~234]{[9]} as the number $e(B)$ of $x \in {\bf n}\setminus B$ for which there is a circuit $C$ such that
$$C \subseteq B \cup \{x\} \quad\mbox{and} \quad x = \max (C).$$
Dually, the {\it internal activity} of $B$ is the number $i(B)$ of $x \in B$ for which there is a cocircuit $D$ such that for $B':= {\bf n} \setminus B$ we have
$$D \subseteq B'\cup \{x\} \quad \mbox{where} \quad x = \max (D).$$
Of course, getting all circuits and cocircuits of ${\cal M}$ is easier said
than done. For graphic ${\cal M}$ a method is given in \cite{[7]}. When ${\cal M}$ is provided as the column matroid of a matrix $A$ over a finite field $K$, we may proceed as in section 3. For $n \leq 15$ and $K = GF(2)$ this straightforward approach suffices.  

As to getting the bases consider the rank 5 polygon matroid ${\cal M}$ on the
edge set ${\bf 8}$ of the graph:
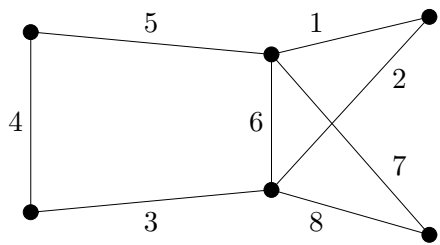
\begin{center}
\setlength{\unitlength}{1mm}
\begin{picture}(60,40)
\put(5,6){\circle*{2}}
\put(5,30){\circle*{2}}
\put(37,9){\circle*{2}}
\put(37,27){\circle*{2}}
\put(58,3){\circle*{2}}
\put(58,32){\circle*{2}}
\drawline (5,6) (5,30)  (37,27) (58,32) (37,9)  (58,3) (37,27) (37,9) (5,6) 
\put(3,18){\makebox(0,0)[c]{4}}
\put(21,30){\makebox(0,0)[b]{5}}
\put(21,6){\makebox(0,0)[t]{3}}
\put(35,18){\makebox(0,0)[c]{6}}
\put(43,6){\makebox(0,0)[t]{8}}
\put(43,30){\makebox(0,0)[b]{1}}
\put(54,25){\makebox(0,0)[t]{2}}
\put(54,11){\makebox(0,0)[b]{7}}
\end{picture}
\end{center}
One reads off that the circuits $C_i$ are (using shorthand notation)
$$126,\quad 678,\quad 1287,\quad 3456,\quad 12345,\quad 34578.$$
Processing them one by one with the $n$-algorithm of \cite{[8]} yields a compact encoding of the family Mod of all $X \subseteq {\bf 8}$ that contain {\it no} set $C_i$. Thus, in our situation, Mod is the family of all independent sets of ${\cal M}$:

\begin{center}
\begin{tabular}{c|c|c|c|c|c|c|c|c|} \hline
& 1 & 2 & 3& 4 & 5& 6 & 7 & 8\\ \hline
$r_1=$ & $n_1$ & $n_1$ & $n_2$ & $n_2$ & $n_2$ & 1 & $n_3$ & $n_3$\\ \hline 
$r_2=$ & $n_1$ & $n_1$ & $n_1$ & $n_1$ & $n_1$ & 0 & $n_2$ & $n_2$\\ \hline
$r_3=$ & $n_1$ & $n_1$ & $n_2$ & $n_2$ & $n_2$ & 0 & 1 & 1 \\ \hline \end{tabular}
\end{center}

Specifically Mod is given as a disjoint union $r_1 \cup r_2 \cup r_3$ where each row $r_i$ consists of certain $0,1$-vectors $X$ of length 8, corresponding to subsets of ${\bf 8}$ in the usual way. The string of symbols (say) $n_2n_2n_2$ in $r_1$ by definition means that only vectors $X$ are allowed that have {\it at least one} $0$ in a position occupied by an $n_2$. Similarly for $n_1n_1$ and $n_3n_3$. Thus $r_1$ contains $(2^2-1)(2^3-1)(2^2-1) = 63$ vectors $X$, one of them say $01101100$. There is a systematic (recursive) procedure to determine the number of sets of fixed cardinality in such a {\it multivalued} row. But here we proceed ad hoc. Namely, the sought bases of being the maximal independent sets,  each base occurs as a maximal member of some row $r_i$. Here, incidentally, each row-maximal member also is a base (since it has cardinality 5). The row-maximal members of $r_i$ are obtained by filling each $n$-bubble to full capacity with $1$'s in all possible ways. For instance, $n_1n_1n_1n_1n_1$ in $r_2$ can be filled with four $1$'s in 5 ways. It follows that $r_2$ features $5 \cdot 2 = 10$ bases. Altogether there are
$$12 + 10 + 6 = 28$$
bases $B$. It is plausible that for each base $B$ the parameters $e(B)$ required in (2) can now be calculated fast. The numbers $i(B)$ can be obtained the same way; after all they are just the numbers $e(B')$ for the dual bases $B'$.

We mention that the bottleneck in the sketched method is getting the circuits and cocircuits (the bases and cobases are then readily found as sketched). How the circuits can be retrieved faster than by browsing the whole cycle space, and how the cocircuits arise {\it from} the circuits (even the chordless ones), is work in progress.

\section{Lists of non-isomorphic regular matroids}
Complete lists with tables of numbers of non-isomorphic regular matroids can 
be found on the web:
\begin{center}
\tt http://www.uni-graz.at/\~{}fripert/html/matroids/matroide\_neu.html
\end{center}
In addition to these tables we also provide complete lists of representatives
of these matroids, and their Tutte polynomials.

As it was previously explained, first we determine lists of non-isomorphic
binary matroids. 
For $1\leq n\leq 15$ and $1\leq k\leq 7$ there are complete lists of
the representatives of the isomorphism classes 
of all (a) loopless, (b) simple loopless, (c) connected simple loopless, and
(d) connected loopless binary matroids of rank $k$ and size $n$.  As the
standard representative of a matroid we choose the lexicographic smallest
representative in its isomorphism class.

By testing these matroids for regularity, we obtain complete lists of
standard representatives of the isomorphism classes of
(a) connected simple loopless and 
(b) connected loopless regular matroids of rank $k$ and size $n$ again for 
$1\leq n\leq 15$ and $1\leq k\leq 7$. 

The lists of matroids of rank $\geq 8$ were computed by changing from a
matroid to its dual. If $\mathcal{M}$ is a binary matroid of size $n$ and rank
$k$, then its dual, $\mathcal{M}^d$ is a binary matroid of size $n$ and rank
$n-k$. Two binary matroids are isomorphic if and only if their duals are
isomorphic. Moreover, a matroid is regular if and only if its dual is regular.
Using this method we determined complete lists of (in general non-standard)
representatives of the isomorphism classes of
(a) connected simple loopless and 
(b) connected loopless regular matroids of rank $k\geq 8$ and size $n\leq 15$.

Additionally, the Tutte polynomials of all connected simple loopless regular 
matroids are determined and can be downloaded in form of Mathematica
Notebooks.


\begin{thebibliography}{10}

\bibitem{ColRead79} C.J. Colbourn and R.C. Read. Orderly algorithms for
generating restricted classes of graphs. {\it Journal of Graph Theory}, 3: 
187--195, 1979.

\bibitem{ColRead79a} C.J. Colbourn and R.C. Read. Orderly algorithms for graph
generation. {\it Int. J. Comput. Math.,} 7: 167--172, 1979.

\bibitem {Grund90} R. Grund. Symmetrieklassen von Abbildungen und die
Konstruktion von diskreten Strukturen. {\it Bayreuth. Math. Schr.}, 31:
19--54, 1990. ISSN 0172-1062.

\bibitem{[6]} G. Haggard, D. Pearce, G. Royle. Computing Tutte polynomials.
{\it ACM Trans. Math. Software},  37, Art. 24, 17pp, 2010.

\bibitem{Kerber99} A. Kerber. {\it Applied Finite Group Actions}, volume 19 of
{\it Algorithms and Combinatorics}. Springer, Berlin, Heidelberg, New York,
1999. ISBN 3-540-65941-2.

\bibitem{[5]}  J.G. Oxley. {\it Matroid Theory}. Oxford Graduate Texts in
Mathematics 3, second edition, 1997.

\bibitem{Read78}
R.C. Read.
\newblock Every one a winner or how to avoid isomorphism search when
  cataloguing combinatorial configurations.
\newblock {\em Ann. Discrete Mathematics}, 2:107 -- 120, 1978.

\bibitem{Symmetrica}
{SYMMETRICA}.
\newblock A program system devoted to representation theory, invariant theory
  and combinatorics of finite symmetric groups and related classes of groups.
  Copyright by ``Lehrstuhl {II} f{\"u}r {M}a\-the\-ma\-tik, Universit{\"a}t
  Bayreuth, 95440 Bayreuth".\\
  \verb+http://www.algorithm.uni-bayreuth.de/en/research/SYMMETRICA/+.

\bibitem{[9]} N. White (ed.) {\it Matroid Applications}. Enc. of Math. and
Appl., Vol.40, Cambridge University Press, 1992.

\bibitem{[7]} M. Wild. Generating all cycles, chordless cycles and Hamiltonian
cycles with the principle of exclusion. {\it Journal of Discrete Algorithms},
6: 93--102, 2008.

\bibitem{[8]} M. Wild. Compactly generating all satisfying truth assignments
of a Horn formula, to appear in {\it Journal on Satisfiability, Boolean
Modeling and Computation}.

\end{thebibliography}
\end{document}